\documentclass[a4paper,12pt]{article}
\usepackage{amsfonts}
\usepackage{amssymb}
\usepackage{amsmath}
\usepackage{amsthm}
\usepackage{hyperref}

\providecommand{\doi}[1]{Eprint \href{http://dx.doi.org/#1}{doi:#1}}
\usepackage{algorithm}
\usepackage{pgf}
\usepackage{tikz}
\usetikzlibrary{arrows}
\usepgflibrary{arrows}
\usetikzlibrary{automata}
\usetikzlibrary{shapes,snakes}
\title{Computing the smallest fixed point
of order-preserving nonexpansive mappings arising in positive stochastic games and static analysis of
programs}
\author{Assal\'e Adj\'e\thanks{LSV, CNRS \& ENS de Cachan, 61, avenue du Pr\'esident Wilson, 
F-94235 Cachan Cedex, France. 
Email: assale.adje@lsv.ens-cachan.fr.
This work was performed when the first author was with the MeASI team of CEA, LIST and with CMAP,
\'Ecole Polytechnique, being supported by a PhD Fellowship of the R\'egion \^Ile-de-France.},
St\'ephane Gaubert\thanks{INRIA Saclay \& Centre de Math\'ematiques Appliqu\'ees, Ecole Polytechnique, 91128 Palaiseau, France.
Email: stephane.gaubert@inria.fr.
}\ \ and
Eric Goubault\thanks{MeASI, CEA LIST, 91191 Gif-sur-Yvette France.
Email: eric.goubault@cea.fr.}
}

\begin{document}
\newtheorem{theorem}{Theorem}[section]
\newtheorem{proposition}{Proposition}[section]
\newtheorem{corollary}{Corollary}[section]
\newtheorem{lemma}{Lemma}[section]
\theoremstyle{definition}
\newtheorem{definition}{Definition}[section]
\newtheorem{def2}{Definition}[section]
\theoremstyle{remark}
\newtheorem{example}{Example}[section]
\newtheorem{remark}{Remark}[section]

\renewcommand{\thefootnote}{\*}
\newcommand{\nn}{\mathbb N}
\newcommand{\rr}{\mathbb R}
\newcommand{\rd}{\rr^d}
\newcommand{\brd}{\overline{\rr}^d}
\newcommand{\ff}{{\mathbf {Fix}}(f)}
\newcommand{\fu}{f'_u}
\newcommand{\pu}{{\mathbf{Fix}_{|\rd_-}(\fu)}}
\newcommand{\cc}{C}
\newcommand{\rh}{{\rho_{\rd_-}(\fu)}}
\newcommand{\vv}{\mathcal V}
\newcommand{\zz}{\mathbb Z}
\def\norm#1{\mbox{$\| #1 \|$}}
\def\evec{\mathbf{e}}
\def\tagp#1{\tag{\mbox{$P_{#1}$}}}
\maketitle

\begin{abstract}
The problem of computing the smallest fixed point of an order-preserving
map arises in the study of zero-sum positive stochastic games. It also
arises in static analysis of programs by abstract interpretation. In this
context, the discount rate may be negative.
We characterize the minimality of a fixed point in terms of the nonlinear
spectral radius of a certain semidifferential. We apply this
characterization to design a policy iteration algorithm, which applies to the case of
finite state and action spaces. The algorithm returns a locally minimal fixed point,
which turns out to be globally minimal when the discount rate is nonnegative.
\end{abstract}


\footnote{
The authors were partly supported by the Arpege  programme of the French National Agency of Research (ANR), project ``ASOPT'', number ANR-08-SEGI-005 and by the Digiteo project DIM08 ``PASO'' number 3389.
}

\section{Introduction}
Zero-sum repeated games  can be studied classically by means
of dynamic programming or Shapley operators.
When the state space is finite, such an operator is a map
$f$ from $\rd$ to $\rd$, where $d$ is the number of states.
Typically, the operator $f$ can be written as:
$$f_i(x)=\min_{a\in A_i}\max_{b\in B_{i,a}} P_i^{a,b}x +r_i^{a,b}$$
Here, $A_i$ represents the set of actions of Player I (Minimizer) in
state $i$, $B_{i,a}$ represents the set of actions of Player II (Maximizer)
in state $i$ when Player I has just played $a$ (the information
of both players is perfect), $r_i^{a,b}$ is an instantaneous
payment from Player I to Player II, and 
$P_i^{a,b}=(P_{ij}^{ab})_j\in\mathbb{R}^d$ is a substochastic vector, giving the transition probabilities
to the next state, as a function of the current state and of
the actions of both players. The difference $1-\sum_j P_{ij}^{a,b}$
gives the probability that the game terminates as a function
of the current state and actions. The operator $f$ will send
$\rd$ to $\rd$ if for instance the instantaneous payments
are bounded.
We may consider the game in which the total payment
is the expectation of the sum of the instantaneous payments of Player I to Player II,
up to the time at which the game terminates. This includes
the discounted case, in which for all $i$, $\sum_j P_{ij}^{a,b}=\alpha<1$,
for some discount factor $\alpha$.
Then, the fixed point of $f$ is unique, and its $i$th-coordinate gives the value
of the game when the initial state is $i$, see~\cite{FV}. In more general
situations~\cite{FV,MS}, the value is known to be the largest (or dually, the smallest) fixed point 
of certain Shapley operators, and it is of interest to compute
this value, a difficulty being that Shapley operators may have
several fixed points.

The same problem appears
in a different context. Static analysis of programs
by abstract interpretation \cite{PCRC} is a technique
to compute automatically invariants of programs, in order to prove
them correct. The fixed point operators arising in static analysis include the
Shapley operators of stochastic games as special cases. However,
the ``discount factor'' may be larger than one, 
which is somehow unfamiliar from the game theoretic point of view.
In this context, the existence of the smallest fixed point
is guaranteed by Tarski-type fixed point arguments,
and this fixed point is generally obtained by 
a monotone iteration (also called Kleene iteration) of the operator
$f$. This method is often slow. Some accelerations based
on ``widening'' and ``narrowing''\cite{NarrWid} are commonly used,
which may lead to a loss of precision, since they only 
yield an upper bound of the minimal fixed point.
Some of the authors introduced
alternative algorithms
based on policy iteration instead \cite{CAV05, ESOP07}, which are often
faster and more accurate. However, the fixed
point that is returned is not always the smallest one.

In the present paper, we refine these policy
iteration algorithms in order to reach the smallest fixed point of $f$
even in degenerate situations. 
Our main result, Theorem~\ref{thm2} below, characterizes the minimality
of a given fixed point in terms of the spectral radius of its semidifferential
map. This is inspired by a result of Akian, Gaubert and Nussbaum~\cite{Fix},
showing that a given fixed point of a semidifferentiable nonexpansive map is unique if and only if 
its semidifferential as $0$ as a unique fixed point (actually, the result of~\cite{Fix} is proved 
in an infinite dimensional setting, using only a mild compactness condition). Theorem~\ref{thm2} 
also shows that when a fixed point is locally minimal (meaning there is no smaller fixed point
in a neighborhood), it is globally minimal. Thus, the present paper
shows that the ideas of localization via semidifferentials developed
in~\cite{Fix} also allow one to address the minimality
issue for a fixed point, instead of the uniqueness issue.

The construction of the present policy iteration
algorithm relies on Theorem~\ref{thm2}, since an eigenvector of the semidifferential
map is used as a descent direction to determine the new policy in degenerate
iterations. 


An alternative
approach to compute the smallest fixed point has been developed by Gawlitza and Seidl~\cite{seidl,DBLP:conf/sas/GawlitzaS10}. 
In a nutshell, Gawlitza and Seidl use an approach dual
to the one of~\cite{CAV05,ESOP07}, iterating in the ``max'' strategy
space instead of the ``min'' strategy space. The advantage
of this  approach
of~\cite{seidl,DBLP:conf/sas/GawlitzaS10} is that it allows one to compute the smallest fixed
point in cases in which the map is expansive in the sup-norm, which
are beyond the scope of the present approach. However, the present
 ``min'' strategy approach has other interests
(in particular, it gives at any step of the algorithm a safe
upper bound of the fixed point, which can be useful in situations
in which the convergence is slow, which do occur in applications).
A detailed comparison
of the two approaches can be found in~\cite{ConvOpt}. 
The questions of computing fixed points of monotone maps, motivated
by verification problems, has also been addressed in~\cite{LEROUX-SUTRE-SAS2007,LEROUX-SUTRE-FSTTCS2007} 
and~\cite{esparza:approximative}.


The present work builds on
the ``operator'' 
(nonexpansive maps) approach to optimal control and zero-sum games, 
which has been developed by several authors, 
see~\cite{neymansurv} for a survey and~\cite{rosenbergsorin,sorin,spectral,vigeral,AGL10,renault,AGGut10} 
for more or less recent developments in this direction.

Finally, we note that an initial version of the present 
work has appeared in the proceedings~\cite{MTNS08},
and that a further version appeared in the Phd thesis of the first author~\cite[Chapter 8]{adjephd}.

\section{Basic notions}
In this paper, we will work in $\rd$ equipped with the sup-norm $\norm{\cdot}$.
We consider the natural partial 
order on $\rd$ defined as: $x\leqslant y$ if for all $i$, $1 \leqslant i \leqslant d$, 
$x_i\leqslant y_i$
where $x_i$ indicates the $i$th coordinate of $x$. We write $x<y$ when $x\leqslant y$ and
there exists a $j$ such that $x_j<y_j$. We denote by $\mathbb{R}_+$ (resp. $\mathbb{R}_-$)
the set of real nonnegative (resp. nonpositive) numbers.
All vectors of $\rd$ such that $f(x)=x$ are called fixed points of the map $f$.
We denote by $\ff$ the set of fixed points of $f$.





When the action spaces are finite, dynamic programming
operators are not differentiable, and they may even
have empty subdifferentials or superdifferentials.
However, following~\cite{Fix}, we analyze their local behavior by means
of a nonlinear analogue of the differential, the
semidifferential (e.g.~\cite{RW}). In order to define it, let
us recall some basic notions concerning cones and homogeneous
maps.





\begin{def2}[Cone, homogeneous map and spectral radius]
A subset $\cc$ of $\rd$ is called a \textit{cone} if for all $\lambda\geqslant 0$
and for all $x\in\cc$, $\lambda x\in\cc$. A cone $\cc$ is said to be \textit{pointed} if 
$\cc\cap-\cc=\{0\}$. A self-map $g$ on a cone $\cc$ is said to be \textit{homogeneous} 
(of degree one), if for all strictly positive real numbers $\lambda$, $g(\lambda x)=\lambda g(x)$. 
We define the \textit{spectral radius} of a homogeneous continuous self-map $g$ on a 
closed convex pointed cone $\cc$, the nonnegative number $\rho_{\cc}(g)$: 
\begin{displaymath}
\rho_{\cc}(g)=\sup\{\lambda\geqslant 0\mid\exists\ x\in\cc\backslash\{0\}, g(x)=\lambda x\}
\end{displaymath}
\end{def2}

Note that closed convex pointed cones are precisely what we need to define spectral
radius of homogeneous continuous self-maps $g$ on $\cc$.  Several notions 
of spectral radius for continuous maps have been defined, we refer the reader
to~\cite{Nus,Eig} for more details. 

A vector $x\in\cc\backslash\{0\}$ such that $g(x)=\lambda x$ is a {\em nonlinear eigenvector} of $g$, and $\lambda$ is 
the associated {\em nonlinear eigenvalue}. 
Remark that, when $g$ maps from $C$ to itself and $C$ is a pointed cone, all the eigenvalues are 
nonnegative. The existence of nonlinear eigenvectors is guaranteed by standard fixed point arguments~\cite{Nus}.  

\begin{proposition}[Positivity of nonlinear spectral radius]
\label{existence_spec}
For any pointed convex cone $\cc$, for any continuous map $g:\cc\mapsto \cc$, the set:
\[
\{\lambda\geqslant 0\mid \exists\, x\in\cc\backslash\{0\},\ g(x)=\lambda x\}
\]
is nonempty.
\end{proposition}


 
We next recall the notion of semidifferential, see~\cite{RW} for more
background.
\begin{definition}[Semidifferential]
Let $u\in\rd$ and $f$ be a self-map on $\rd$. We say that 
$f$ is semidifferentiable at $u$ if there exists a homogeneous continuous
map $g$ on $\rd$ and a neighborhood $\vv$ of 0 such that for all $h\in\vv$:
\begin{displaymath}
f(u+h)=f(u)+g(h)+o(\|h\|)
\end{displaymath}
We call $g$ the semidifferential of $f$ at $u$ and we denote it $\fu$.
\end{definition}
Note that $f$ semidifferentiable at $u$ implies that $f$ is continuous at $u$.  
If $f$ is semidifferentiable at $u$, we have for all $t>0$ and for $h$ small enough:
$$f(u+th)=f(u)+t\fu(h)+o(t\norm{h})$$
This implies that $$\fu(h)=\lim_{t\to 0^+}\frac{f(u+th)-f(u)}{t}$$ and
the semidifferential coincides with the directional derivative of $f$ at $u$
in direction $h$ (on the positive side). The latter limit characterization implies 
that the semidifferential is unique. The following result shows that semidifferentiability
requires the latter limit to be uniform in the direction $h$.

\begin{proposition}[{see~\cite[Theorem 7.21]{RW}}]
\label{propRW}
Let $f$ be a self-map on $\rd$. Let $u$ be in $\rd$.
The map $f$ is semidifferentiable at $u$ if and only if
for all vectors $h$, the following limit exists:
\begin{displaymath}
\lim_{\substack{t\to 0^+\\h'\to h}}\dfrac{f(u+th')-f(u)}{t} \enspace 
\end{displaymath}
\end{proposition}

In this paper, we also treat the special case of piecewise affine maps.

\begin{definition}[Piecewise affine map]
\label{affmor}
Let $f$ be a self-map on $\rd$. We say that $f$ is piecewise affine if 
for all $j\in\{1,\cdots,d\}$ there exists finite sets $A_j$ and $\{B_a\}_{a\in A_j}$ and
a family $\{g_{a,b}\}$ of affine maps such that:
\[
f_j=\min_{a\in A_j}\max_{b\in B_a} g_{a,b}
\]
\end{definition}

It is shown in~\cite{Ovc,CONEDUALITY} 
that the set of piecewise affine maps that we define is the same
as the set of functions $f$ for which there exists a family of convex closed sets with nonempty interior
which covers $\rd$ and such that the restriction of $f$ on each element 
of this family is affine.

\begin{proposition}
\label{calsem}
Let $f$ be a piecewise affine self-map on $\rd$. Let $u$ be in $\rd$. Then $f$ is semidifferentiable
at $u\in\rd$. For $j\in\{1,\ldots,d\}$, we set $\overline{A_j}=\{a\in A_j\mid f_j(u)=\max_{b\in B_a} g_{a,b}(u)\}$ and 
$\overline{B_a}=\{b\in B_a\mid g_{a,b}(u)=\max_{b\in B_a} g_{a,b}(u)\}$, then:
\[
(\fu)_j=\min_{a\in\overline{A_j}}\max_{b\in\overline{B_a}} \nabla g_{a,b}\cdot
\]
\end{proposition}

This assertion may be deduced by applying the rule of the ``differentiation'' of a max 
see Exercise~10.27 of \cite{RW}. Alternatively Proposition~\ref{calsem} can be recovered 
from Proposition~\ref{directiondesc} (established independently) below.


\section{Main Results}

Let $f$ be a continuous map from $\rd$ to $\rd$ with a nonempty set of fixed points.
For $u\in\rd$, we denote the set of nonpositive fixed points of $\fu$ by $\pu$. 

We recall a simple case where an order-preserving self-map $f$ on $\rd$ 
has a smallest fixed point in $\rd$.

\begin{proposition}[Existence of smallest fixed point]
\label{existence}
Let $f$ be a self-map on $\rd$. Suppose that $f$ is order-preserving.
Assume that the set $\{x\in\rd\mid f(x)\leqslant x\}$ is bounded from below. 
Then $f$ has a smallest fixed point $u$ in $\rd$ and $u$ satisfies:
\[
u=\inf\{x\in\rd\mid f(x)\leq x\} \enspace .
\] 
\end{proposition}

\begin{proof}
Let $\Gamma$ be the set $\{x\in\rd\mid f(x)\leqslant x\}$. The fact that $\Gamma$ is bounded
from below implies that the infimum of $\Gamma$ exists, we denote it by $u$. 
Since $f$ is order-preserving, we have $f(u)\leqslant u$ and $f(u)\in\Gamma$. 
Hence we get $u\leq f(u)$ and so $u=f(u)$. The set $\Gamma$ contains all the fixed points of $f$ 
and we conclude that $u$ is the smallest fixed of $f$ in $\rd$. 
\end{proof}

The existence of the smallest fixed point for $f$ does not imply that 
the set $\{x\in\rd\mid f(x)\leqslant x\}$ is bounded from below. Let us take for 
example, the function on $\rr$, $g:x\mapsto 2x$. The function $g$ has a unique fixed 
point in $\rr$ but the $\{x\in\rr\mid g(x)\leqslant x\}=\rr_-$. 
Proposition~\ref{proptarskinonexp} below will describe a situation where the existence 
of a smallest fixed point for a map $f$ implies that the set 
$\{x\in\rd\mid f(x)\leqslant x\}$ is bounded from below.  

\begin{definition}[Locally minimal fixed point]
Let $u\in\ff$. We say that $u$ is a locally minimal fixed point if there is 
a neighborhood $\vv$ of $u$ such that for all $v\in\vv\cap\ff$, 
$v\leqslant u \implies v=u$. 
\end{definition}

\begin{theorem}
\label{thm1}
Let $u\in\ff$. Assume that $f$ is semidifferentiable at $u$.
Consider the following statements:
\begin{enumerate}
 \item [1.] $u$ is a locally minimal fixed point.
 \item [2.] $\pu=\{0\}$.
 \item [3.] $\rh<1$. 
\end{enumerate}

Then $3\ \implies\ 2\ \implies 1$.
\end{theorem}

\begin{proof}
Point $3$ implies point $2$ indeed.
In order to 
show that $2$ implies $1$, assume that $u$ is not a locally minimal
fixed point. Then, there exists a sequence $h_n$ of 
nonzero vectors in $\rd_-$ tending to the zero vector such that $u+h_n$
is a fixed point of $f$. After replacing $h_n$ by a subsequence,
we may assume that $y_n:=\norm{h_n}^{-1}h_n$ has a limit $y$.
Then, $\|y\|=1$ and $y\in \rd_-$.
Writing $u+h_n=u+\norm{h_n}y_n$,
and using Proposition~\ref{propRW}, we conclude that $y=\fu(y)$,
showing that $\fu$ has a nonzero fixed point in $\rd_-$.
\end{proof}

\begin{example}
\label{localex}
Consider the self-map $f$ of $\rr^2$ defined as follows:
\[
f
\left(
\begin{array}{c}
x\\
y
\end{array}
\right)
=
\left(
\begin{array}{c}
\max(\min(x,y),0)+2\max(\min(x+1,0),-2)\\
\max(\max(x,y),0)+2\max(\min(y+1,0),-2)
\end{array}
\right)
\]

The fixed point set of $f$ is shown at Figure~\ref{locale}. It consists of the union 
of six isolated points, two half lines and a cone.
Let us deduce the local minimality of the fixed point $(0,0)$ by using 
Theorem~\ref{thm1}.

\begin{figure}[ht!]
\begin{center}
\resizebox{9cm}{!}{\includegraphics{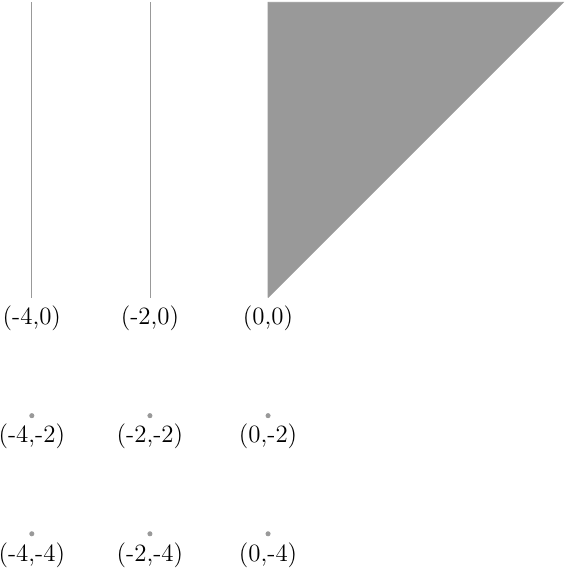}}
\end{center}
\caption{Fixed point set of the function $f$ of Example~\ref{localex}}
\label{locale}
\end{figure}

Using Proposition~\ref{calsem}, the semidifferential at $(0,0)$ is the function:
\[
h\mapsto f'_{(0,0)}(h_1,h_2)=
(\max(\min(h_1,h_2),0),\max(\max(h_1,h_2),0)),
\]
 hence, $f'_{(0,0)}$ vanishes on the nonpositive 
cone of $\rr^2$ and so $f'_{(0,0)}$ does not have any nonzero negative fixed point. 
\end{example}

We just gave sufficient condition to the local minimality. Now, we give a sufficient 
condition to ensure that a fixed point is not the smallest one.

\begin{proposition}
\label{thm1b}
Suppose that $f$ is an order-preserving self-map and that the set $\{x\in\rd\mid f(x)\leqslant x\}$ 
is bounded from below. Let $u$ be in $\ff$ and assume that $f$ is semidifferentiable at $u$. 
If $\rh>1$ then $u$ is not the smallest fixed point.
\end{proposition}

\begin{proof}
 By definition of the spectral radius, $\rh>1$ implies that, there exists $\mu>1$ and $h\in \rd_-\backslash\{0\}$ such that,
$\fu(h)=\mu h$. Let us write $\epsilon:=\displaystyle{\min_{h_i<0}-\norm{h}^{-1}(\mu-1)h_i}>0$. There exists $t_0>0$ such that,
$0<t<t_0$ implies that $f_i(u+th)\leqslant u_i+t\mu h_i+t\epsilon\norm{h}$, for all $i$ such that $h_i<0$. Hence, for all $i$ 
such that $h_i<0$, $f_i(u+th)\leqslant u_i+t\mu h_i-t\norm{h}^{-1}(\mu-1)h_i\norm{h}=u_i+th_i<u_i$. Since $f$ is order-preserving,
$f(u+th)\leqslant u$, for all $t\geqslant 0$, so $f_i(u+th)\leqslant u_i+th_i$, for all $i$ such that $h_i=0$.

We conclude that, $f(u+th)<u+th$. The sequence $(f^{k}(u+th))_{k\in\mathbb{N}}$ is nonincreasing and all its terms belong 
to the set $\{x\in\rd\mid f(x)\leqslant x\}$ which is bounded from below, by continuity of $f$, it converges to a fixed point strictly smaller than $u$.
\end{proof}

\begin{example}
We come back to Example~\ref{localex}. The map $f$ is clearly bounded from below
since $f(x,y)\geqslant (-4,-4)$ for all $(x,y)\in\rr^2$. We consider the fixed point $(0,-2)$. 
At this point, the semidifferential $f'_{(0,-2)}$ is the function:
\[
h\mapsto (0,\max(h_1,0)+2h_2),
\]
thus, $\rho_{\rr_-^2}(f'_{(0,-2)})=2$. By Theorem~\ref{thm1b}, $(0,-2)$ is not the smallest fixed point of 
$f$.
\end{example}

In Theorem~\ref{thm1} and Proposition~\ref{thm1b}, there are no restrictive conditions on the map. 
We next consider the special case of piecewise affine maps.

\begin{proposition}
\label{directiondesc}
Let $f$ be a piecewise affine self-map on $\rd$. There is neighborhood $\vv$ of $u$ such that, for all $u+h\in\vv$:
\[
f(u+h)=f(u)+\fu(h)\enspace .
\]
\end{proposition}

\begin{proof}
Note that it is equivalent to show the equality for all coordinates of $f$ and 
we fix a coordinate $j\in\{1,\ldots,d\}$. We set, for all $a\in A_j$, $g_a(x)=\max_{b\in B_a} g_{a,b}(x)$.
We denote by $\overline{B_a}$ the set of elements $b\in B_a$ such that $g_a(u)=g_{a,b}(u)$. 
There exists a neighborhood $\vv_a$ of $u$ such that: $g_a(u+h)=\max_{b\in\overline{B_a}} g_{a,b}(u+h)$, 
for all $h$ such that $u+h\in\vv_a$. Since $g_{a,b}$ is affine, we have: 
$g_{a,b}(u+h)=g_{a,b}(u)+\nabla g_{a,b}\cdot h$.
It follows that, for all $\bar b\in\overline{B_a}$:
\begin{equation}
g_a(u+h)=g_{a,\bar b}(u)+\max_{b\in\overline{B_a}}\nabla g_{a,b}\cdot h
\label{eq1}
\end{equation}
Let us denote $\overline{A_j}$ the set of elements $a\in A_j$ such that $f_j(u)=g_a(x)$. 
There exists a neighborhood of $u$, $\vv\subseteq\bigcap_a\vv_a$, such that:
$f_j(u+h)=\min_{a\in\overline{A_j}} g_a(u+h)$ if $u+h\in\vv$.
Applying (\ref{eq1}), we get:
$f_j(u+h)=f_j(u)+\min_{a\in\overline{A_j}}\max_{b\in\overline{B_a}}\nabla g_{a,b}\cdot h$ if $u+h\in\vv$.{}
Clearly $h\mapsto \min_{a\in\overline{A_j}}\max_{b\in\overline{B_a}}\nabla g_{a,b}\cdot h$ is continuous and homogeneous
thus it must coincide with the semidifferential of $f_j$ at $u$. Note that 
Proposition~\ref{calsem} follows from the latter fact.
\end{proof}

From Proposition~\ref{directiondesc}, we deduce the following corollary which characterizes 
the locally minimal fixed point of piecewise affine maps.

\begin{corollary}
\label{desc}
Let $f$ be a piecewise affine self-map on $\rd$ and let $u\in\ff$, then
$u$ is a locally minimal fixed point if and only $\pu=\{0\}$.
\end{corollary}

The second part of the previous proposition is the basis
of a ``descent'' algorithm given in the section \ref{algor}. 
If a fixed point $u$ is not locally minimal, then there exists a strictly
negative fixed point $h$ for $\fu$ which may be thought of as a 
{\em descent direction} such that $u+h$ is a fixed point of $f$. 

\begin{example}
We consider the same example as before, Example~\ref{localex}. Now, we look at a fixed 
point of the form $z=(x,y)$ such that $y\geqslant x>0$. The map $f$ of Example~\ref{localex} is piecewise
affine. Following Proposition~\ref{calsem}, at $z$ the semidifferential $f'_z$ is the function:
\[
h\mapsto\left\{
\begin{array}{cc}
(h_1,h_2)&\text{ if } y>x\\
(\min(h_1,h_2),\max(h_1,h_2)) &\text{ if } y=x
\end{array}
\right.
\] 

In these two cases, $f'_z$ admits a negative fixed point hence $z$ is not a locally fixed point.
The isolated points of Figure~\ref{locale} and the vector $(0,0)$ do not have any negative fixed point,
we recover that these fixed points are locally minimal. 
\end{example}

In order to pass from local minimality to global minimality,
we shall need the following nonexpansiveness condition.
\begin{definition}[Nonexpansive map]
Let $f$ be a self-map on $\rd$: $f$ is nonexpansive (with respect to the
sup-norm) if for all $x,y\in\rd$, $\norm{f(x)-f(y)}\leqslant \norm{x-y}$.
\end{definition}

Nonexpansiveness is automatically satisfied by Shapley operators since they 
are order-preserving and additively homogeneous.

\begin{proposition}
\label{ray}
Let $f$ be a nonexpansive self-map on $\rd$ and let $u\in\rd$. Assume that $f$ is semidifferentiable 
at $u$. Then $\pu=\{0\}$ if and only if $\rh<1$.
\end{proposition}
\begin{proof}
It suffices to show that $\rh\leqslant 1$.
Firstly, $\fu$ is nonexpansive as a
pointwise limit of nonexpansive maps.
Assume that $\rh>1$, there exists
$\mu>1$ and $v\in\rd_-$ such that $\fu(v)=\mu v$,
then $\norm{\fu(v)-\fu(0)}=\mu\norm{v}>\norm{v-0}$ which 
contradicts the nonexpansiveness of $\fu$.
\end{proof}

The main result of this paper is the following theorem, which
will allow us to check the global minimality of a fixed point.

\begin{theorem}
\label{thm2}
\label{fond}
Let $f$ be an order-preserving nonexpansive self-map on $\rd$.
Let $u$ be a fixed point of $f$. 
Then, $u$ is locally minimal if and only if it is 
the smallest fixed point of $f$.
If in addition $f$ is piecewise affine,
the following assertions are equivalent:
\begin{enumerate}
 \item [1.] $u$ is the smallest fixed point of $f$.
 \item [2.] $\pu=\{0\}$.
 \item [3.] $\rh<1$.
\end{enumerate}
\end{theorem}

The proof of this theorem relies on the existence of an order-preserving 
and nonexpansive retract on the fixed point set of $f$. This idea
was already used in \cite{CAV05}.
The existence of nonexpansive retracts on the fixed point set is
a classical topic in the theory of nonexpansive mappings,
see Nussbaum \cite{Nus}. In the present finite dimensional case,
the result of the next lemma is an elementary one.

\begin{lemma}
\label{previous}
Let $f$ be a nonexpansive order-preserving self-map on $\rd$. Let $u$ be in $\ff$.
Then there is a nonlinear order-preserving and nonexpansive map $P$ such that
$P(\rd)=\ff$ and $\mathbf{Fix }(P)=\ff$.
\end{lemma}
\begin{proof}
Following the idea of \cite{Per}, we shall construct the map 
$P$ as follow:
$P(x)=\lim_{k\to +\infty} f^k(y)$ where 
$y=\limsup_{l\to +\infty} f^l(x)$.
Since $f$ is nonexpansive and $u\in\ff$, $(f^k(x))_{k\in\nn}$ is bounded for all $x\in\rd$. We can now write,
for all $x\in\rd$, $Q(x)=\limsup_{k\to +\infty}f^k(x)$.
Moreover, given $k\geqslant 0$, we have, for all $m\geqslant k$,
$\sup_{n\geqslant k}f^n(x)\geqslant f^m(x)$ and since $f$ is order-preserving, for all $m\geqslant k$,
$f(\sup_{n\geqslant k}f^n(x))\geqslant f(f^m(x))$ 
so, $f(\sup_{n\geqslant k}f^n(x))\geqslant \sup_{m\geqslant k}f(f^m(x))$,
we conclude, by taking the limit when $k$ tends to $+\infty$ and using
the continuity of $f$, $f(Q(x))\geqslant Q(x)$, so $(f^l(Q(x))_{k\in\nn}$ is a nondecreasing sequence.
Moreover, since $f$ is nonexpansive, the limit $P(x)=\lim_{l\to +\infty}f^l(Q(x))$ is finite. Observe that the map 
$P$ is order-preserving and nonexpansive since it is the pointwise limit of order-preserving and nonexpansive maps.
Furthermore, $f(P(x))=f(\lim f^l(Q(x))=\lim f^{l+1}Q(x)=P(x)$ so $P(\rd)\subseteq\ff$.
Moreover, it is easy to see that $P$ fixes every fixed point of $f$. It follows that
$P$ is a projector.
\end{proof}

\begin{proof}[Proof of Theorem~\ref{fond}]
Suppose that $u$ is a locally minimal fixed point but not the smallest fixed point.
Then, there is a fixed point $v$ such that $\inf(v,u)<u$. 
For all scalars $t\geqslant 0$, define 
$\omega_t:= \inf(v+t,u)$.
Let us take $P$ as in Lemma~\ref{previous}.
Since $P$ is nonexpansive in the sup-norm,
$\norm{P(v+t)-P(v)}\leqslant t$ for all $t\geqslant 0$ and so
$P(v+t)\leqslant P(v)+t$.
Using the monotonicity of $P$, we deduce
that $P(\omega_t)\leqslant \inf(P(v+t),P(u))\leqslant \inf(P(v)+t,P(u))
=\omega_t$. 
Let $t_0=\inf\{t\geqslant 0\mid\omega_t=u\}$. Then, for $0<t<t_0$,
$P(\omega_t)$ is a fixed point of $f$, which is such that
$P(\omega_t)<u$. Since $P$ is continuous, $P(\omega_t)$ tends to $P(\omega_{t_0})=P(u)=u$ as $t$ tends to $t_0^-$,
which contradicts the local
minimality of $u$. Hence a contradiction with $u$ not being the smallest fixed
point. Define by $1'$ the property that $u$ is a locally minimal fixed point.
We just showed $1\Leftrightarrow 1'$. By Theorem~\ref{thm1}, we get
$3\implies 2\implies 1'$. By Corollary~\ref{desc}, $1'\implies 2$ and 
by Proposition~\ref{ray} we get $2\implies 3$.  
\end{proof}

\section{A policy iteration algorithm to compute the smallest fixed point}
\label{algor}
The previous results justify the following policy iteration algorithm
which returns the smallest fixed point of a nonexpansive order-preserving 
piecewise affine map. 
Assume that $f$ is a map from $\rd$ to $\rd$, every coordinate
of which is given by 
\begin{equation}
\label{selec}
f_j(x)=\inf_{a\in A_j} f^{a}(x)
\end{equation}
where $A_j$ is finite and every $f_a$ is a supremum of order-preserving nonexpansive affine maps.
A strategy $\pi$ is a map from $\{1,\cdots,d\}$ to 
$\mathcal A=\bigcup_{1\leqslant j\leqslant d} A_j$ such that $\pi(j)\in A_j$ for all $j$.
We define $f^{\pi}=(f^{\pi(1)},\cdots,f^{\pi(j)},\cdots,f^{\pi(d)})$. We assume that every map 
$f^{\pi}$ has a smallest fixed point.
The idea of the algorithm is to use a descent direction to select the new strategy when a non-minimal 
fixed point is reached .
The algorithm, needs two oracles. Oracle 1 returns the smallest
fixed point in $\rd$ of a map $f^{\pi}$. Oracle 2 checks whether the restriction of $\fu$ to
the convex cone $\rd_-$ has a spectral radius equal to 1, and if this is the case, 
returns a vector $h\in\rd_-\backslash\{0\}$ such that $\fu(h)=h$. 
We discuss below the implementation of these oracles for subclasses of maps.

\begin{algorithm}
\caption{Computing the smallest fixed point by Policy Iteration} \label{policymod}
\textsl{Input}: An order-preserving nonexpansive map $f$ in the form \eqref{selec}.\\
\textsl{Output}: The smallest fixed point of $f$ in $\rd$.\\
\textbf{Init}:   Select a strategy $\pi^0$, $k=0$.\\
\textbf{Value Determination} $\mathbf{(D_k)}$:
Call Oracle 1 to compute the smallest fixed point $u^k$ of $f^{\pi^k}$ in $\rd$.\\
\textbf{Policy Improvement} $\mathbf{(I_k^1)}$:
If $f(u^k)<u^k$, take $\pi^{k+1}$ such that $f(u^k)=f^{\pi^{k+1}}(u^k)$ and go 
to Step $\mathbf{(D_{k+1})}$.\\
\textbf{Policy Improvement} $\mathbf{(I_k^2)}$:
If $f(u^k)=u^k$, call Oracle 2 to compute $\alpha_k:=\rho(f_{u^k}^{\prime})$.

$\bullet$ If $\alpha_k<1$, return $u^k$, which is the smallest fixed point of $f$. 

$\bullet$ If $\alpha_k=1$, take $h\in\rd_-\backslash\{0\}$ such that $f_{u^k}^{\prime}(h)=h$. 
Define $\pi^{k+1}(j)$ as an optimal action $a$ in $(f_{u^k}^{\prime})_j(h)=\min_{a\in \overline{A_j}}
(f^{a})^{\prime}_{u^k}(h)$ where $\overline{A_j}=\{a\in A_{j}\mid f^{a}(u^k)=f_{j}(u^k)\}$. 
Then go to $\mathbf{(D_{k+1})}$.
\end{algorithm}

Since the number of strategies is finite, we shall see that the existence of a smallest fixed point in $\rd$ 
for every policy suffices to show that Algorithm~\ref{policymod} terminates. We start by proving 
that the smallest fixed point of a nonexpansive and order-preserving 
self-map on $\rd$ can be computed using an optimization problem.    

\begin{proposition}\label{proptarskinonexp}
Let $g$ be an order-preserving nonexpansive self-map on $\rd$.
Assume that $g$ has a smallest fixed point in $\rd$ denoted by $u$. 
\begin{enumerate}
\item The vector $u$ is the smallest element of $\{x\in\rd\mid g(x)\leqslant x\}$;
\item The vector $u$ is the unique optimal solution of the minimization problem:
      $\min\{\sum_{1\leqslant i\leqslant d} x_i\mid x\in\rd,\ g(x)\leqslant x\}$. 
\end{enumerate}
\end{proposition}

\begin{proof}
The vector $u$ satisfies $g(u)=u$ and a fortiori $g(u)\leqslant u$. Let $x$ be 
any vector satisfying $g(x)\leqslant x$. Since $g$ is order-preserving, nonexpansive and has a fixed point, 
$(g^k(x))_{k\in\nn}$ is a bounded nonincreasing sequence. Denoting by $y$ the limit of $(g^k(x))_{k\in\nn}$, we get 
$y\leqslant x$ and $y$ is a finite fixed point of $g$ so $y\geqslant u$, the first assertion is thus proved.
We conclude next that $\sum_{1\leqslant i\leqslant d} u_i\leqslant 
\sum_{1\leqslant i\leqslant d} y_i\leqslant \sum_{1\leqslant i\leqslant d} x_i$. 
Since this holds for all feasible $x$, it follows that $u$ is an optimal solution. If $x$ is an arbitrary optimal solution, 
we must have $\sum_{1\leqslant i\leqslant d} u_i=\sum_{1\leqslant i\leqslant d} x_i$ and since $u\leqslant x$, 
it follows that $u=x$ and the second assertion is proved.
\end{proof}

Now, we prove that Algorithm~\ref{policymod} terminates when policies have a smallest fixed point.

\begin{theorem}[Termination of Algorithm~\ref{policymod}]
Let $f$ be an order-preserving nonexpansive map in the form~\eqref{selec}.
If every policy $f^{\pi}$ has a smallest fixed point in $\rd$, then Algorithm~\ref{policymod} terminates.
\end{theorem}

\begin{proof}
To show that Algorithm~\ref{policymod} terminates, it suffices to check that the sequence of 
produced points $u^0,u^1,\cdots$ is strictly decreasing, because the corresponding policies
must be distinct and the number of policies is finite. Let $k$ be an integer. Let be $u^k$ be a vector of $\rd$. 
Suppose that an improvement of type $I_k^1$ arises: $f(u^k)<u^k$. There exists a policy $f^{\pi^{k+1}}$ such that 
$f^{\pi^{k+1}}(u^k)=f(u^k)$. Let us denote $u^{k+1}$ the smallest fixed point of $f^{\pi^{k+1}}$.
Since $u^k$ belongs to the set $\{x\in\rd\mid f^{\pi^{k+1}}(x)\leqslant x\}$, we conclude from the first assertion 
of Proposition~\ref{proptarskinonexp} that $u^{k+1}\leqslant u^k$. Furthermore, since 
$u^k$ is not a fixed point of $f^{\pi^{k+1}}$, $u^{k+1}<u^k$.
Suppose that an improvement of type $I_k^2$ arises: $u^k$ is not the smallest fixed point of $f$.
Then, since $f$ is piecewise affine and nonexpansive, by Proposition~\ref{desc} and Theorem~\ref{fond}, 
there exists $h\in\rd_-\backslash \{0\}$, $f^{\pi^{k+1}}(u^k+th)=u^k+th$ for $t>0$ small enough. 
It follows that $u^{k+1}\leqslant u^k+th<u^k$.
\end{proof}

Now, we discuss the implementation of the oracles.
Using Proposition~\ref{proptarskinonexp} assertion 2, the following corollary identifies 
a situation where Oracle 1 can be implemented by solving a linear programming problem.

\begin{corollary}\label{prop-lp}
Let $g$ be an order-preserving nonexpansive map that is the
supremum of finitely many affine maps.
Assume, that $g$ has a smallest fixed point $u$ in $\rd$. Then $u$ is the unique 
optimal solution of the linear program:
\[
\min\{\sum_{1\leqslant i\leqslant d} x_i\mid x\in\rd,\ g(x)\leqslant x\}.
\]
\end{corollary}

The implementation of Oracle 2 raises the issue of computing the spectral radius.
Let $g$ be an order-preserving, homogeneous and continuous self-map of $\rd_-$. 
It is known that:
\begin{align}
\rho_{\rd_-}(g)&=\inf_{x\in\operatorname{int}(\rd_-)}\sup_{1\leqslant i\leqslant d}\frac{g_i(x)}{x_i}
\label{spec1}\\
\rho_{\rd_-}(g)&=\sup_{x\in\mathbf{\rd_-}}\limsup_{k\to +\infty}\norm{g^k(x)}^{\frac{1}{k}}
\label{spec2}
\end{align}
The first equality, which is a generalization of the Collatz-Wielandt property
in Perron-Frobenius theory, follows from a result of Nussbaum \cite{Cvx} Theorem~3.1.
The second characterization is shown by Mallet-Paret and Nussbaum in \cite{Eig} under more general assumptions.
We deduce, for every vector $x\in\operatorname{int}(\rd_-)$:
\begin{equation}
\rho_{\rd_-}(g)\leqslant\left(\sup_{1\leqslant i\leqslant d}\frac{g_i^k(x)}{x_i}\right)^{\frac{1}{k}}.
\label{spec0}
\end{equation}
Moreover, the latter upper bound converges to $\rho_{\rd_-}(g)$ as $k$ tends to infinity.
This yields an obvious method to check whether $\rho_{\rd_-}(g)<1$, which consists in 
computing the upper bound in \eqref{spec0} for successive values of $k$ as long as the 
upper bound is not smaller than 1. This algorithm will not stop when $\rho_{\rd_-}(g)=1$.
However, we next describe simple situations where this idea leads to a terminating algorithm.

First, we can adapt Proposition~\ref{prop-lp} to compute a descent direction from a linear program 
when the semidifferential coincides with a supremum of finitely many affine maps. Since the 
semidifferential is a homogeneous map we can add box constraints to the linear program of
Proposition~\ref{prop-lp}. From now, we denote $\evec$ the vector the coordinates of which are equal to
$-1$.

\begin{proposition}\label{prop-lphomogene}
Let $g$ be a homogeneous order-preserving nonexpansive self-map.
Assume, that the restriction of $g$ on the nonpositive cone $\rd_-$ is the supremum of finitely many affine maps.
Let $\alpha$ be a nonnegative real number. 
The unique optimal solution of the linear program:\\
\begin{equation}
\tagp{\alpha}
\label{lphomogene}
\min\left\{\sum_{1\leqslant i\leqslant d} x_i\mid x\in\rd,\ g(x)\leqslant x,\ \alpha\evec\leqslant x\leqslant 0\right\}
\end{equation}
is the smallest fixed point of $g$ in the set $\{x\in\rd\mid \alpha\evec\leqslant x\leqslant 0\}$.

Furthermore, if there exists $\alpha>0$ such that the optimal value of~\eqref{lphomogene} is equal to 0 
then the null vector is the unique fixed point of $g$ in $\rd_-$.
\end{proposition}

\begin{proof}
Since $g$ is nonexpansive and homogeneous, $\norm{g(\alpha\evec)}\leq \alpha$ and thus 
$g(\alpha\evec)\geqslant \alpha\evec$. Moreover, $g$ is order-preserving, we conclude that $g$ is a self-map 
on $\{x\in\rd\mid \alpha\evec\leqslant x\leqslant 0\}$ which is a complete lattice. From Tarski's theorem, $g$ 
has a smallest fixed point $u$ such that $\alpha\evec\leqslant u\leqslant 0$ and $u$ is the smallest 
element of the feasible set of~\eqref{lphomogene}. Using the same argument as in assertion 2 of 
Proposition~\ref{proptarskinonexp}, we conclude that $u$ is the unique optimal solution of~\eqref{lphomogene}.

Now, suppose there exists $\alpha>0$ such that the optimal value of~\eqref{lphomogene} is 0 and assume that $g$
has a nonzero fixed point $u$ in $\rd_-$. Since $g$ is homogeneous, $\alpha u/\norm{u}$ is a feasible 
point of~\eqref{lphomogene}, the sum of coordinates of this vector is strictly smaller than 0 which 
contradicts the fact that 0 is the optimal value of~\eqref{lphomogene}.   
\end{proof}

The second simple case concerns the homogeneous min-max functions.

\begin{definition}
\label{minmax}
We call a homogeneous min-max function of the variables $h_1,\ldots,h_d$ a term in the grammar:
$X\mapsto\min(X,X),\max(X,X),h_1,\cdots,h_d,0$. 
\end{definition}

For instance, the term $\min(h_1,\max(h_2,h_3,0))$ is produced by this grammar. More generally, 
we shall say that a map from $\rd$ to $\rd$ is a homogeneous min-max self-map
if its coordinates are of the form of Definition~\ref{minmax}. This definition is inspired by the min-max
functions considered by Gunawardena \cite{MMax} and Olsder \cite{Eigminmax}. The terms of this form 
comprise the semidifferentials of the min-max functions considered there. For simple classes of programs, 
like the one we shall consider in the next section, the semidifferential at any fixed point
turns out to be a homogeneous min-max function. In this case, the spectral radius can be 
computed efficiently by using to the following integrity argument $g(\zz^d)\subseteq\zz^d$.

\begin{proposition}
\label{caleff}
Let $g$ be a homogeneous min-max self-map on $\rd$. Then
$\rho_{\rd_-}(g)\in\{0,1\}$. Moreover $\rho_{\rd_-}(g)=0$ if and only if $\lim_{k\to +\infty} g^k(\evec)=0$, 
and the latter limit is reached in at most $d$ steps.
\end{proposition}

\begin{proof}
Since $g$ is nonexpansive and $g(0)=0$, we have $g(\evec)\geqslant \evec$.
We deduce, by monotonicity of $g$, that $(g^k(\evec))_{k\in\nn}$ is a nondecreasing sequence bounded 
from above by 0. Moreover, $g$ preserves the set of integer vectors.
So this sequence converges in at most $d$ steps to some vector $b\in\zz_-^d$.
Let us suppose that $b=0$.
For all $y\in\rd_-$, there exists $t\geqslant 0$ such that, $0\geqslant y\geqslant t\evec$. Since $g$ is homogeneous 
and order-preserving, $0\geqslant g^k(y)\geqslant tg^k(\evec)=0$ for all $k\geqslant d$ which implies
that $\rho_{\rd_-}(g)=0$. If $b<0$ we have $g(b)=b\neq 0$, and so $\rho_{\rd_-}(g)\geqslant 1$. Finally, since
$g$ is nonexpansive, we also have $\rho_{\rd_-}(g)\leqslant 1$.
\end{proof}

\section{Applications}
\label{example}

\subsection{Applications to positive stochastic games}

In this subsection, we consider applications to game theory, particularly in zero-sum positive
stochastic games. The term "positive" means that all the rewards are nonnegative.  
We are interested in the expectation of the infinite sum of rewards at each date without discounting.
Filar and Vrieze~\cite[Theorem 4.4.3]{FV} prove that the value of a positive stochastic game is the smallest nonnegative 
fixed point of a Shapley operator. They also develop a value iteration~\cite[Theorem 4.4.4]{FV} to compute 
the value of this class of games. Although value iteration is generally an effective method, one can construct 
examples in which policy iteration is much faster. This occurs when the contraction rate of value iteration 
is close to one, or when the initial estimate of the value function is too far away from the true value function. 
We give such an example in Example~\ref{slowexample}. In consequence, we propose to use Algorithm~\ref{policymod} to compute 
the value of zero-sum positive stochastic game with perfect information.


\begin{example}
\label{slowexample}
Let us consider the positive game in perfect information with three states represented by Figure~\ref{slowexamplefig}. 
The real $\epsilon$ is chosen in $(0,1/2]$. In the game the two players can stop the game at each state they can play. 
The diamond node (state 1) is a stochastic node, the circle node (state 2) is controlled by $MIN$ and the triangle node 
(state 3) is controlled by $MAX$. The thick edges represent the deterministic transitions, whereas dotted ones represent 
stochastic transitions. The double octagons are the stopping options for each player and the numbers inside are stopping payoffs. 
The numbers on the edges are either immediate payoffs if the nodes are controlled by players, or the probability 
to reach the state at the endpoint of the edge. 

\begin{figure}[ht!]
\begin{center}
\includegraphics{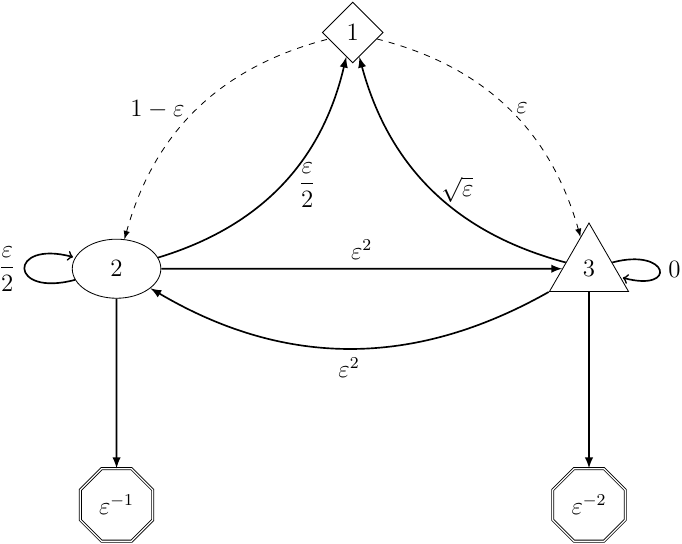}
\end{center}
\label{slowexamplefig}
\caption{The game of Example~\ref{slowexample}}
\end{figure}

The dynamic operator $f_{\varepsilon}$ is thus defined as the following function:
\begin{displaymath} 
f_{\varepsilon}(x)=
f_{\varepsilon}\left(\begin{array}{c}
          x_1\\
          x_2\\
          x_3
    \end{array}\right)
           =
\left(\begin{array}{ccc}
(1-\varepsilon)x_2+\varepsilon x_3\\
     \min(x_1+\dfrac{\varepsilon}{2},x_2+\dfrac{\varepsilon}{2},x_3+\varepsilon^2,\varepsilon^{-1})\\
     \max(x_1+\sqrt{\varepsilon},x_2+\varepsilon^2,x_3,\varepsilon^{-2})\\ 
        \end{array}\right) 
\end{displaymath}

Since $\varepsilon\in (0,1/2]$, the value of the game is 
$V_{\varepsilon}=(2\varepsilon^{-1}-1,\varepsilon^{-1},\varepsilon^{-2})$. Indeed, $V_{\varepsilon}$
is a fixed point of $f_{\varepsilon}$ and the semidifferential at $V_{\epsilon}$ is the following map:
$(h_1,h_2,h_3)\mapsto ((1-\varepsilon)h_2+\varepsilon h_3,0,\max(h_3,0))$ which only admits $(0,0,0)$
as nonpositive fixed point and using Theorem~\ref{fond}, we conclude that $V_{\epsilon}$ is 
the smallest (positive) fixed point of $f_{\varepsilon}$ and thus the value of the game.
An optimal policy is given by selecting the stopping options for each player. Since 
we chose these stopping options as initial policy (this is our heuristic choice for the next examples),
policy iteration algorithm does not need any further iteration to converge. 

Table~\ref{iterationnbtable} shows the numbers of iterations needed by the value iteration 
(VI) for different values of $\varepsilon$. The value iteration is initialized 
with the null vector $(0,0,0)$ and stops when the $\ell_{\infty}$ norm of 
the difference between two iterates are smaller than $10^{-6}$.

\begin{table}[ht!]
\centering
\begin{tabular}{|c|c|}
\hline
$\varepsilon$ & \# iterations for VI \\
\hline
0.5 & 11\\
\hline
0.25 & 35\\
\hline
0.1 & 203\\
\hline 
0.01 & 200,003 \\
\hline
0.001 & 2,000,003 \\
\hline
\end{tabular}
\caption{Numbers of iterations for Value Iteration for different values of $\varepsilon$} 
\label{iterationnbtable}
\end{table}
\end{example}

\subsubsection{A detailed example of a positive game in perfect information} 

\begin{figure}[ht!]
\begin{center}
\resizebox{14cm}{!}{\includegraphics{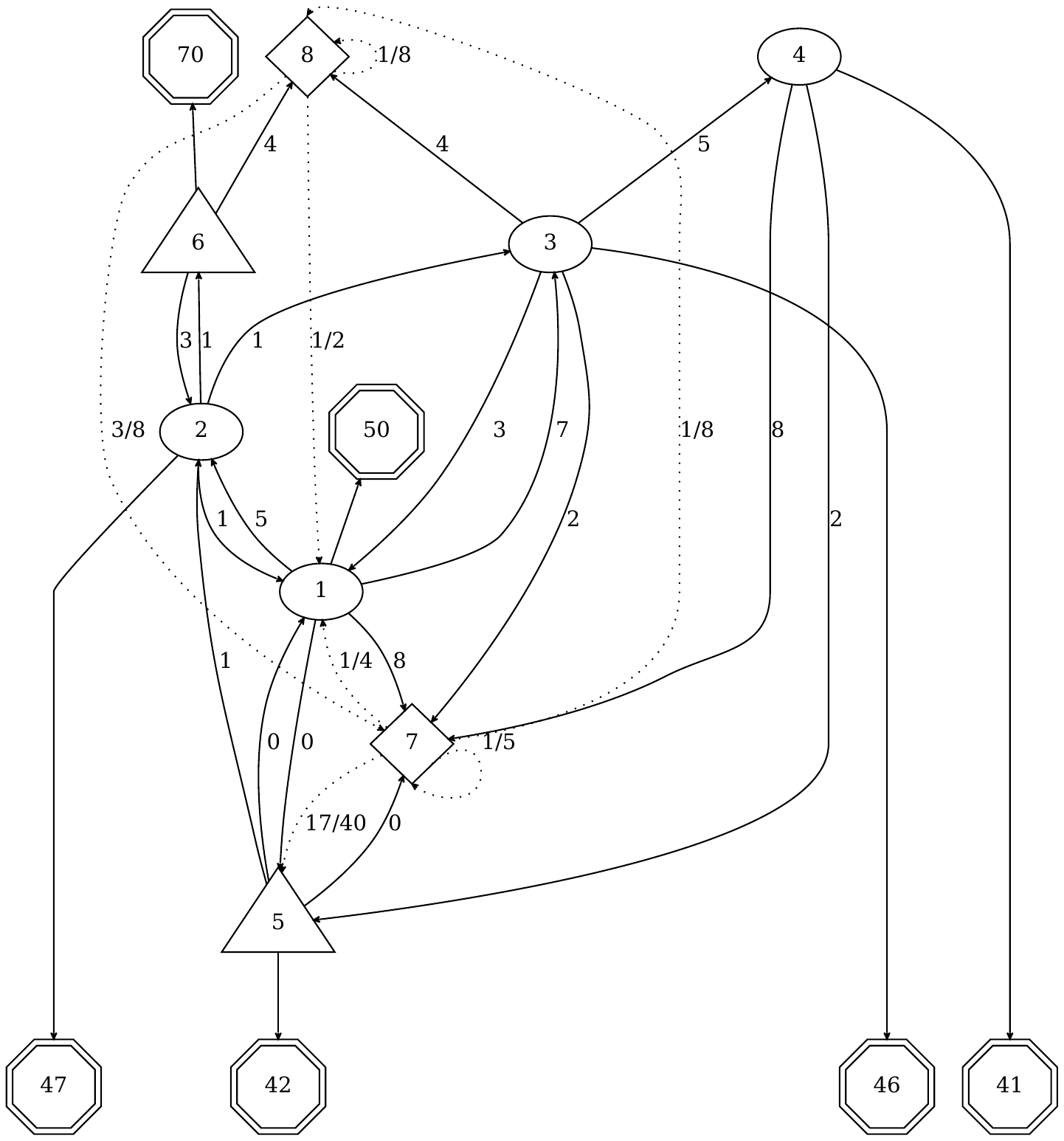}}
\end{center}
\caption{The graph associated of the dynamic operator}
\label{examplegame}
\end{figure}

As in Example~\ref{slowexample}, we consider a positive stochastic games with stopping options.
Again the space of states are decomposed in three kinds of states: the states controlled by $MIN$, the ones 
controlled by $MAX$ and the stochastic states. The graph depicted by Figure~\ref{examplegame} represents the 
game. The circle states are controlled by $MIN$, the triangle ones by $MAX$ and the diamond ones 
are the stochastic states. The other nodes (represented by double octagons) are the stopping options and the numbers 
inside represent the payoff. For the nodes controlled by the players, the thick edges are the different actions 
for the players and the numbers on them are the immediate payoff. The dotted edges are the different transitions 
from a stochastic state to other states. 

The Shapley operator for this game is given by the following order-preserving piecewise affine nonexpansive map:

\begin{displaymath} 
f(x)=
f\left(\begin{array}{c}
          x_1\\
          x_2\\
          x_3\\
          x_4\\
          x_5\\
          x_6\\
          x_7\\
          x_8
\end{array}\right)
           =
\left(\begin{array}{ccc}
     \min(x_2+5,x_3+7,x_5,x_7+8,50)\\
     \min(x_1+1,x_3+1,x_6+1,47)\\
     \min(x_1+3,x_4+5,x_7+2,x_8+4,46)\\
     \min(x_5+2,x_7+8,41)\\
     \max(x_1,x_2+1,x_7,42)\\ 
     \max(x_2+3,x_8+4,70)\\
     (1/4)x_1+(17/40)x_5+(1/5)x_7+(1/8)x_8\\
     (1/2)x_1+(3/8)x_7+(1/8)x_8       
        \end{array}\right) 
\end{displaymath}

We can compute the value of the game by value iteration, it takes 42 iterations.  
Here, we apply Algorithm~\ref{policymod} to compute the smallest nonnegative fixed point of $f$
which determines the value of the game. To do this we can use a map $\tilde{f}$ which coincides with 
the Shapley operator $f$ on the nonnegative cone and the fixed point of which are nonnegative. To construct 
$\tilde{f}$, it suffices to replace for the states controlled by $MIN$ and the stochastic states,
the linear forms $a\cdot x+b_i$ by $\max(a\cdot x+b_i,0)$. Now, we initialize (Step $\mathbf{(D_0)}$) 
Algorithm~\ref{policymod} by taking the stopping options. We get the following map $\tilde{f}^0$:
\begin{displaymath} 
\tilde{f}^{\pi^{0}}(x)
           =
\left(\begin{array}{ccc}
     50\\
     47\\
     46\\
     41\\
     \max(x_1,x_2+1,x_7,42)\\ 
     \max(x_2+3,x_8+4,70)\\
     \max((1/4)x_1+(17/40)x_5+(1/5)x_7+(1/8)x_8,0)\\
     \max((1/2)x_1+(3/8)x_7+(1/8)x_8,0)       
 \end{array}\right) 
\end{displaymath}

By a linear program (Proposition~\ref{prop-lp}), we compute the smallest fixed 
point of $\tilde{f}^{\pi^{0}}$, we find the vector $\bar{x}=(50,47,46,41,50,70,50,50)^{\intercal}$.
As required by step $\mathbf{(I_0^1)}$, we have to check whether $\bar{x}$ is 
fixed point of $\tilde{f}$, this is the case.  
Now we enter in step $\mathbf{(I_1^1)}$, to check whether $\bar{x}$ is the smallest fixed point of $\tilde{f}$, 
we compute the semidifferential of $\tilde{f}$ at $\bar{x}$ in the direction $h$
using Proposition~\ref{calsem}, we obtain the map $\tilde{f}'_{\bar{x}}$
defined as follows:

\[
\left(\begin{array}{ccc}
     \min(h_5,0)\\
     \min(h_3,0)\\
     \min(h_4,0)\\
          0\\
     \max(h_1,h_7)\\ 
          0\\
     (1/4)h_1+(17/40)h_5+(1/5)h_7+(1/8)h_8\\
     (1/2)h_1+(3/8)h_7+(1/8)h_8       
 \end{array}\right)
\] 

On the nonpositive cone, the components of $\tilde{f}'_{\bar{x}}$ are linear forms 
except for the sixth component which is a maximum of linear forms. Then, we can use 
Proposition~\ref{prop-lphomogene} as Oracle~2 to determine whether $\bar{x}$ is the smallest fixed point of $\tilde{f}$.
The unique optimal solution of linear program 
$\min\{\sum_j h_j\mid \tilde{f}'_{\bar{x}}(h)\leqslant h, h_i\in [-1,0],\ \forall i\}$ is 
$\bar{h}=(-1,0,0,0,-1,0,-1,-1)^{\intercal}$. 
We use this vector as descent direction. The vector $\bar{h}$ gives us a new policy which is:

\begin{displaymath} 
\tilde{f}^{\pi^{1}}(x)=
\left(\begin{array}{ccc}
     x_5\\
     47\\
     46\\
     41\\
     \max(x_1,x_2+1,x_7,42)\\ 
     \max(x_2+3,x_8+4,70)\\
     (1/4)x_1+(17/40)x_5+(1/5)x_7+(1/8)x_8\\
     (1/2)x_1+(3/8)x_7+(1/8)x_8       
        \end{array}\right) 
\end{displaymath}

and we go to the step $\mathbf{(D_1)}$. By a linear program (Proposition~\ref{prop-lp}), we compute the smallest fixed 
point of $\tilde{f}^{\pi^{1}}$, we find the vector $x^*=(48,47,46,41,48,$ \- $70,48,48)^{\intercal}$. The vector $x^*$ is 
also a fixed point of $\tilde{f}$. We enter in the step $\mathbf{(I_1^1)}$. We compute the semidifferential of $\tilde{f}$ at 
$x^*$ in the direction $h$ using Proposition~\ref{calsem}, we obtain the map $\tilde{f}'_{x^*}$
defined as follows:

\[
\hspace{1cm}
\left(\begin{array}{ccc}
     \min(h_5,0)\\
     \min(h_3,0)\\
     \min(h_4,0)\\
          0\\
     \max(h_1,h_2,h_7)\\ 
          0\\
     (1/4)h_1+(17/40)h_5+(1/5)h_7+(1/8)h_8\\
     (1/2)h_1+(3/8)h_7+(1/8)h_8       
 \end{array}\right)
\] 

To check whether $x^*$ is the smallest fixed point of $\tilde{f}$, we can proceed by computing 
the upper bound of Equation~\eqref{spec0} where $g$ is the semidifferential of $\tilde{f}$ 
at $x^*$. Taking the vector $e$ the coordinates of which are equal to $-1$, we find
that $\max_{1\leqslant i\leqslant 8} \dfrac{g_i^6(e)}{e_i}=\dfrac{109}{320}<1$ and we conclude that
$x^*$ is the smallest fixed point of $\tilde{f}$ and so the value of the game 
depicted by Figure~\ref{examplegame}.
 
\begin{remark}
We could also use the same linear program~\eqref{lphomogene} with $\alpha=1$. 
In the case of the semidifferential of $\tilde{f}$ at $x^*$, the optimal value would be equal 
to zero and the unique optimal solution would be the null vector. 
\end{remark}

\subsection{Applications to static analysis of programs}

We next illustrate our results on an example from
static analysis. We take a simple but interesting program with nested
loops (Figure~\ref{program}). 
{From} this program, we create semantic equations
on the lattice of intervals \cite{PCRC} (Figure~\ref{equations}) that
describe the outer approximations of the sets of values that program
variables can take, for all possible executions. 
For instance, at control point 4, the value of variable $y$ can come
either from point 3 or point 5 (hence the union operator), 
as long as the condition $y\geqslant 5$ is
satisfied (hence the intersection operator).
An interval $I$ is written
as $[-i^-,i^+]$ in order to get fixed point equations of order-preserving maps
in $i^i$ and $i^+$.
The equations we derive on bounds are
order-preserving piecewise affine maps, to which we can apply our methods.

\begin{figure}[ht!]
\begin{minipage}{0,4\textwidth}
\begin{tabular}{l c}
\verb|int x,int y,|     &\\ 
\verb|x=[0,2];y=[10,15]|&//1\\
\verb|while (x<=y) {     |&//2\\
\verb|  x=x+1;     |&//3\\
\verb|  while (5<=y) {|&//4\\
\verb|    y=y-1;|&//5\\
\verb|  }    |&//6\\
\verb|}          |&//7
\end{tabular}
\caption{A simple C program}
\label{program}
\end{minipage}
\begin{minipage}{0,4\textwidth}
\begin{eqnarray*}
(x_1,y_1)&=&([0,2],[10,15])\\
x_2&=&(x_1\cup x_6)\cap[-\infty,(y_1\cup y_6)^+]\\
y_2&=&(y_1\cup y_6)\cap[(x_1\cup x_6)^-,+\infty]\\
(x_3,y_3)&=&(x_2+[1,1],y_2)\\
(x_4,y_4)&=&(x_3,(y_3\cup y_5)\cap [5,+\infty])\\
(x_5,y_5)&=&(x_4,y_4+[-1,-1])\\
(x_6,y_6)&=&(x_5,(y_3\cup y_5)\cap [-\infty,4])\\
x_7&=&(x_1\cup x_6)\cap[(y_1\cup y_6)^-+1,+\infty]\\
y_7&=&(y_1\cup y_6)\cap[-\infty,(x_1\cup x_6)^+-1]
\end{eqnarray*}
\caption{Its abstract semantic equations in intervals}
\label{equations}
\end{minipage}
\end{figure}

The order-preserving nonexpansive 
piecewise affine map $f$ for the bounds of these intervals is:
\begin{displaymath} 
f\left(\begin{array}{c}
          x\\
          y
       \end{array}\right)
           =
f\left(\begin{array}{c}
          x_2^-\\
          x_2^+\\
          x_7^-\\
          x_7^+\\
          y_2^-\\
          y_2^+\\
          y_4^-\\
          y_4^+\\
          y_6^-\\
          y_6^+\\
          y_7^-\\
          y_7^+
\end{array}\right)
           =
\left(\begin{array}{ccc}
                0&\vee&(x_2^--1)\\
     2\vee (x_2^++1)&\wedge& \underline{15\vee y_6^+}\\
     0\vee (x_2^--1)&\wedge&\underline{(-10\vee y_6^-)-1}\\
                 0&\vee&(x_2^++1)\\
     \underline{0\vee (x_2^--1)}&\wedge& -10\vee y_6^-\\ 
                15&\vee& y_6^+\\
            y_2^-\vee (y_4^-+1)&\wedge&\underline{-5}\\
             y_2^+&\vee& y_4^+-1\\
                y_2^-&\vee& y_4^-+1\\
             y_2^+\vee (y_4^+-1)&\wedge&\underline{4}\\
             -10&\vee&y_6^-\\
         15\vee y_6^+&\wedge&\underline{(2\vee (x_2^++1))-1} 
        \end{array}\right) 
\end{displaymath}

In the equations for the intervals $x_2,y_2,y_4,y_6,x_7$ and $y_7$,
an intersection appears, which gives a $\min$ ($\wedge$) in the corresponding
coordinate of $f$. Choosing a policy is the same as replacing
every minimum of terms by one of the terms, which yields
a simpler ``minimum-free'' nonlinear map, which can
be interpreted as the dynamic programming operator
of a one-player problem.

We next illustrate Algorithm~\ref{policymod}. 
The underlined terms in the expression of $f$ indicate the initial policy
$\pi^0$, for instance, the fifth coordinate of $f^{\pi^0}$ 
is $0\vee (x_2^--1)$.
We first compute the smallest fixed point of this policy
(Step $\mathbf{(D_0)}$).
This may be done by linear programming (Proposition~\ref{prop-lp}),
or, in this special case, by a reduction to a shortest
path problem.
We find $(\bar x,\bar y)=(0,15,-1,16,0,15,-5,15,0,4,0,15)^{\intercal}$.
The first Policy Improvement step, $\mathbf{(I^1_0)}$, requires
to check whether this is a fixed point of $f$. This
turns out to be the case. To determine whether $(\bar x,\bar y)$ 
is actually the smallest fixed point, we enter in the second
policy improvement step, $\mathbf{(I^2_0)}$.
We calculate the semidifferential at $(\bar x,\bar y)$
in the direction $(\delta x,\delta y)$, using Proposition~\ref{calsem}:

$$\begin{array}{rcl} 
f'_{(\bar x,\bar y)}
   \left(          \delta\bar x, 
           \delta\bar y
       \right)^{\intercal}
       &      = &
\left(
                    0,
                    0,
                 \delta\bar{y}_6^-,
                 \delta\bar{x}_2^+,
              0\wedge\delta\bar {y}_6^-, 
                   0,
                   0,
             \delta\bar{y}_2^+,
           \delta\bar{y}_2^-,
                   0,
             \delta\bar{y}_6^-,
         0\wedge\delta\bar{x}_2^+
        \right)^{\intercal}
\end{array}$$
The method of Proposition~\ref{caleff}, which we use as Oracle~2, 
yields in three steps a fixed point for
$f'_{(\bar x,\bar y)}$ that we denote by $h=(0,0,-1,0,-1,0,0,0,-1,0,-1,0)^{\intercal}$.
This fixed point determines the new policy $\pi^1$, which corresponds
to the map $g:=f^{\pi^1}$ given by:
\begin{displaymath} 
g\left(\begin{array}{c}
          x\\
          y
       \end{array}\right)
           =
g\left(\begin{array}{c}
          x_2^-\\
          x_2^+\\
          x_7^-\\
          x_7^+\\
          y_2^-\\
          y_2^+\\
          y_4^-\\
          y_4^+\\
          y_6^-\\
          y_6^+\\
          y_7^-\\
          y_7^+
\end{array}\right)
           =
\left(\begin{array}{ccc}
                0&\vee& x_2^--1\\
                15&\vee& y_6^+\\
              -11&\vee& y_6^--1\\
                 0&\vee&x_2^++1\\
               -10&\vee& y_6^-\\ 
                15&\vee& y_6^+\\
                   &-5&\\
             y_2^+&\vee& y_4^+-1\\
              y_2^-&\vee& y_4^-+1\\
                    &4&\\
             -10&\vee&y_6^-\\
         1&\vee&x_2^+ 
        \end{array}\right) 
\end{displaymath}

We are now in a new value determination step, $(D_1)$.
We find the fixed point $(\tilde u,\tilde v)=(0,15,-5,16,-4,15,-5,15,-4,4,-4,15)^{\intercal}$ of $g$, which is also a fixed point of $f$
(Step $\mathbf{(I^1_1)}$). In Step $\mathbf{(I^2_1)}$, using Proposition~\ref{calsem}, we compute
the semidifferential of $f$ at $(\tilde u,\tilde v)$, which is given
by:
$$\begin{array}{rcl}
f'_{(\tilde u,\tilde v)}
   \left(
          \delta\tilde u,
           \delta\tilde v
       \right)^{\intercal}
      &      =&
\left(
                    0,
                    0,
             \delta\tilde{v}_6^-,
                 \delta\tilde{u}_2^+,
              \delta\tilde{v}_6^-, 
                   0,
                   0,
             \delta\tilde{v}_2^+,
           \delta\tilde{v}_2^-\vee\delta\tilde{v}_4^-,
                   0,
             \delta\tilde{v}_6^-,
         0\wedge\delta\tilde{u}_2^+ 
        \right)^{\intercal}
\end{array}$$

Calling again Oracle~2, we get $\rho_{\rd_-}(f'_{(\tilde u,\tilde v)})=0$,
and so, the algorithm stops: $(\tilde u,\tilde v)$ is the smallest fixed point of $f$.

\bibliographystyle{alpha}
\bibliography{arxivbib}

\newcommand{\etalchar}[1]{$^{#1}$}
\begin{thebibliography}{CGG{\etalchar{+}}05}

\bibitem[Adj11]{adjephd}
A.~Adj\'e.
\newblock {\em Optimisation et jeux appliqu\'es \`a l'analyse statique de
  programme par interpr\'etation abstraite}.
\newblock Phd thesis, \'Ecole Polytechnique, April 2011.

\bibitem[AG03]{spectral}
M.~Akian and S.~Gaubert.
\newblock Spectral theorem for convex monotone homogeneous maps, and ergodic
  control.
\newblock {\em Nonlinear Analysis. Theory, Methods \& Applications},
  52(2):637--679, 2003.
\newblock \arxiv{math.SP/0110108}, \doi{10.1016/S0362-546X(02)00170-0}.

\bibitem[AGG08]{MTNS08}
A.~Adj\'e, S.~Gaubert, and E.~Goubault.
\newblock Computing the smallest fixed point of nonexpansive mappings arising
  in game theory and static analysis of programs.
\newblock In {\em Proceedings of the Eighteenth International Symposium on
  Mathematical Theory of Networks and Systems (MTNS2008)}, Blacksburg,
  Virginia, July 2008.
\newblock \arxiv{0806.1160}.

\bibitem[AGG11]{AGGut10}
M.~Akian, S.~Gaubert, and A.~Guterman.
\newblock Tropical polyhedra are equivalent to mean payoff games.
\newblock To appear in International of Algebra and Computation,
  \arxiv{0912.2462}, \doi{10.1142/S0218196711006674}, 2011.

\bibitem[AGL11]{AGL10}
M.~Akian, S.~Gaubert, and B.~Lemmens.
\newblock Stability and convergence in discrete convex monotone dynamical
  systems.
\newblock {\em Journal of Fixed Point Theory and Applications}, 9(2):295--325,
  2011.
\newblock \doi{10.1007/s11784-011-0052-1}, \arxiv{1003.5346}.

\bibitem[AGN]{Fix}
M.~Akian, S.~Gaubert, and R.D. Nussbaum.
\newblock Uniqueness of fixed point of nonexpansive semidifferentiable maps.
\newblock {\em preprint}.
\newblock \arxiv{1201.1536}.

\bibitem[AT07]{CONEDUALITY}
C.~D. Aliprantis and R.~Tourky.
\newblock {\em Cone and duality}.
\newblock AMS, 2007.

\bibitem[CC77]{PCRC}
P.~Cousot and R.~Cousot.
\newblock Abstract interpretation: A unified lattice model for static analysis
  of programs by construction of approximations of fixed points.
\newblock {\em Principles of Programming Languages 4}, pages 238--252, 1977.

\bibitem[CC92]{NarrWid}
P.~Cousot and R.~Cousot.
\newblock Comparing the galois connection and widening/narrowing approaches to
  abstract interpretation, invited paper.
\newblock In {\em Proceedings of the International Workshop Programming
  Language Implementation and Logic Programming}, Leuven, Belgium, 13--17
  August 1992, Lecture Notes in Computer Science 631, pages 269--295.
  Springer-Verlag, 1992.

\bibitem[CGG{\etalchar{+}}05]{CAV05}
A.~Costan, S.~Gaubert, E.~Goubault, M.~Martel, and S.~Putot.
\newblock A policy iteration algorithm for computing fixed points in static
  analysis of programs.
\newblock In {\em Proceedings of Computer Aided Verification 2005}, Lecture
  Notes in Computer Science 3576. Springer, august 2005.

\bibitem[EGKS08]{esparza:approximative}
J.~Esparza, T.~Gawlitza, S.~Kiefer, and H.~Seidl.
\newblock Approximative methods for monotone systems of min-max-polynomial
  equations.
\newblock In {\em Proceedings of the 35th international colloquium on Automata,
  Languages and Programming (ICALP'08), Part I}, pages 698--710, 2008.

\bibitem[FV97]{FV}
J.~Filar and K.~Vrieze.
\newblock {\em Competitive Markov Decision Processes}.
\newblock Springer-Verlag, 1997.

\bibitem[GG04]{Per}
S.~Gaubert and J.~Gunawardena.
\newblock The {P}erron-{F}robenius theorem for homogeneous monotone functions.
\newblock {\em Transactions of AMS}, 356(12):4931--4950, 2004.

\bibitem[GGTZ07]{ESOP07}
S.~Gaubert, E.~Goubault, A.~Taly, and S.~Zennou.
\newblock Static analysis by policy iteration on relational domains.
\newblock In {\em Proceedings of European Symposium Of Programming 2007},
  Lecture Notes in Computer Science 4421, pages 237--252. Springer, 2007.

\bibitem[GS07]{seidl}
T.~Gawlitza and H.~Seidl.
\newblock Precise fixpoint computation through strategy iteration.
\newblock In {\em Proceedings of the Sixteenth European Symposium on
  Programming (ESOP'07)}, Lecture Notes in Computer Science 4421, pages
  300--315. Springer, 2007.

\bibitem[GS10]{DBLP:conf/sas/GawlitzaS10}
T.~Gawlitza and H.~Seidl.
\newblock Computing relaxed abstract semantics w.r.t. quadratic zones
  precisely.
\newblock In R.~Cousot and M.~Martel, editors, {\em SAS}, volume 6337 of {\em
  Lecture Notes in Computer Science}, pages 271--286. Springer, 2010.

\bibitem[GSA{\etalchar{+}}12]{ConvOpt}
T.~Gawlitza, H.~Seidl, A.~Adj\'e, S.~Gaubert, and E.~Goubault.
\newblock Abstract interpretation meets convex optimization.
\newblock {\em Journal of Symbolic Computation}, 47(12):1416 -- 1446, 2012.
\newblock International Workshop on Invariant Generation.

\bibitem[Gun94]{MMax}
J.~Gunawardena.
\newblock Min-max functions.
\newblock {\em Discrete Event Dynamic Systems}, 4:377--406, 1994.

\bibitem[LS07a]{LEROUX-SUTRE-SAS2007}
J.~Leroux and G.~Sutre.
\newblock Accelerated data-flow analysis.
\newblock In {\em Static Analysis, 14th International Symposium, SAS 2007,
  Kongens Lyngby, Denmark, August 22-24, 2007, Proceedings}, volume 4634 of
  {\em Lecture Notes in Computer Science}, pages 184--199. Springer, 2007.

\bibitem[LS07b]{LEROUX-SUTRE-FSTTCS2007}
J.~Leroux and G.~Sutre.
\newblock Acceleration in convex data-flow analysis.
\newblock In {\em Foundations of Software Technology and Theoretical Computer
  Science, 27th International Conference, FSTTCS 2007, New Delhi, India,
  December 12-14, 2007, Proceedings}, volume 4855 of {\em Lecture Notes in
  Computer Science}, pages 520--531. Springer, 2007.

\bibitem[MPN02]{Eig}
J.~Mallet-Paret and R.D. Nussbaum.
\newblock Eigenvalues for a class of homogeneous cone maps arising from
  max-plus operators.
\newblock {\em Discrete and Continuous Dynamical Systems}, 8(3):519--562, July
  2002.

\bibitem[MS97]{MS}
A.D. Maitra and W.D. Sudderth.
\newblock Discrete gambling and stochastic games.
\newblock {\em Journal of the Royal Statistical Society. Series A (Statistics
  in Society)}, 160(2):376--377, 1997.

\bibitem[Ney03]{neymansurv}
A.~Neyman.
\newblock Stochastic games and nonexpansive maps.
\newblock In {\em Stochastic games and applications (Stony Brook, NY, 1999)},
  volume 570 of {\em NATO Sci. Ser. C Math. Phys. Sci.}, pages 397--415. Kluwer
  Acad. Publ., Dordrecht, 2003.

\bibitem[Nus86]{Cvx}
R.D. Nussbaum.
\newblock Convexity and log convexity for the spectral radius.
\newblock {\em Linear Algebra And Its Applications}, 73:59--122, 1986.

\bibitem[Nus88]{Nus}
R.D. Nussbaum.
\newblock Hilbert's projective metric and iterated nonlinear maps.
\newblock {\em Memoirs of the AMS}, 75(391), August 1988.

\bibitem[Ols91]{Eigminmax}
G.J. Olsder.
\newblock Eigenvalues of dynamical max-min systems.
\newblock {\em Discrete Event Dynamic Systems}, 1:177--207, 1991.

\bibitem[Ovc02]{Ovc}
S.~Ovchinnikov.
\newblock Max-min representation of piecewise linear functions.
\newblock {\em Contributions to Algebra and Geometry}, 43(1):297--302, 2002.

\bibitem[QR11]{renault}
M.~Quincampoix and J.~Renault.
\newblock On the existence of a limit value in some non expansive optimal
  control problems.
\newblock {\em SIAM Journal on Control and Optimization}, 49:2118--2132,
  October 2011.

\bibitem[RS01]{rosenbergsorin}
D.~Rosenberg and S.~Sorin.
\newblock An operator approach to zero-sum repeated games.
\newblock {\em Israel J. Math.}, 121(1):221--246, 2001.

\bibitem[RW98]{RW}
R.T. Rockafellar and R.J-B. Wets.
\newblock {\em Variational Analysis}.
\newblock Springer, 1998.

\bibitem[Sor02]{sorin}
S.~Sorin.
\newblock {\em A first course on Zero-Sum Repeated Games}.
\newblock Springer, 2002.

\bibitem[Vig10]{vigeral}
G.~Vigeral.
\newblock Evolution equations in discrete and continuous time for nonexpansive
  operators in banach spaces.
\newblock {\em ESAIM: Control, Optimisation and Calculus of Variations},
  16:809--832, 2010.

\end{thebibliography}
\end{document}